\documentclass[10pt]{article}

\usepackage{amsthm,amsmath,amssymb,amsbsy}

\usepackage{graphicx}

\usepackage[colorlinks=true,citecolor=blue,linkcolor=blue,urlcolor=blue]{hyperref}

\theoremstyle{plain}
\newtheorem{theorem}{Theorem}[section]
\newtheorem{lemma}[theorem]{Lemma}

\newtheorem{proposition}[theorem]{Proposition}

\theoremstyle{definition}
\newtheorem{definition}[theorem]{Definition}
\newtheorem{example}[theorem]{Example}

\theoremstyle{remark}
\newtheorem{remark}[theorem]{Remark}

\title{\bf Fundamental groups of 3-dimensional small covers}

\author{Vladimir Gruji\'c\\
\small University of Belgrade, Serbia\\[-0.8ex]
\small Faculty of Mathematics\\[-0.8ex]
\small\tt vgrujic@matf.bg.ac.rs}

\date{
}

\begin{document}

\maketitle

\begin{abstract}
Small covers arising from 3-dimensional simple polytopes are an interesting class of 3-manifolds. The fundamental group is a rigid invariant for wide classes of 3-manifolds, particularly for orientable Haken manifolds, which include orientable small covers. By using Morse-theoretic approach we give a procedure to get an explicit, balanced presentation of the fundamental group of a closed, orientable 3-dimensional simple cover with minimal number of generators. Beside that the minimal Heegaard splitting is determined by this presentation.


\end{abstract}

\section{Introduction}

A small cover $M\xrightarrow{\pi}P$ over a simple $n$-dimensional polytope $P$ is a closed $n$-dimensional manifold $M$ with a locally standard $\mathbb{Z}_2^n$-action such that the orbit space $M/\mathbb{Z}_2^n$ is diffeomorphic to $P$ as a manifold with corners. The examples of small covers include the fixed point sets of complex conjugation on smooth, projective toric varieties and more general on quasitoric manifolds, which are their complex analogue. The main topological properties of small covers are derived in the seminal paper of Davis and Januszkiewicz \cite{DJ}, among others $\mathbb{Z}_2$-cohomology ring is calculated in terms of the Stanley-Reisner ring of the underlying polytope $P$.

The 3-dimensional small covers are an interesting class of 3-manifolds whose topology and geometry attracted attention in recent works. In \cite{BEMPP} is proven the cohomological rigidity of $3$-dimensional small covers over Pogorelov polytopes, a right-angled simple polytopes in Lobachevsky space. Erokhovets built an explicit decomposition of 3-dimensional small covers into geometric parts \cite{E}.

The fundamental group $\pi_1(M)$ of a small cover $M\xrightarrow{\pi}P$ of an arbitrary dimension is presented in \cite{WY}. This presentation, based on the canonical cubical decomposition of a small cover is far to be an optimal. The minimal number of generators of $\pi_1(M)$ is $m-n$, where $m$ is the number of facets of the polytope $P$ and $n$ is its dimension, while the Wu and Yu presentation has $2^{n-1}(m-n)$ generators. They gave a simple combinatorial condition on $\pi_1$-injectivity of proper face submanifolds of a simple cover. In the dimension 3 this condition is satisfied by all small covers over flag polytopes. In the unpublished arXiv version of the paper \cite{WY} is proved that every small covers are Haken manifolds. The classical result of Waldhausen \cite{W} states that the fundamental group determines the diffeomorphism type of an orientable Haken manifold. Therefore orientable 3-dimensional small covers are distinguished up to diffeomorphism by their fundamental groups.

In this paper we describe a procedure for finding an optimal, with the minimal number of generators, presentation of the fundamental group $\pi_1(M)$ of an orientable 3-dimensional small cover $M\xrightarrow{\pi}P$, which is a balanced as well. Our method is based on the Morse-theoretic approach. A Morse function on the small cover $M$ is constructed from a generic linear functional on the underlying polytope $P$ composed with the orbit map. The corresponding handle decomposition of $M$ leads to proposed presentation of the fundamental group.

In section 2 we recall some basic definitions and properties of 3-dimensional small covers. In section 3 we review some basic facts about Haken manifolds and the Heegaard splittings. We show how the Wu and Yu presentation of the fundamental group $\pi_1(M)$ of a small cover $M\xrightarrow{\pi}P$ in the dimension 3 is easily read off the canonical Heegaard splitting associated to the characteristic pair of $M$. In section 4 we describe our method for finding an optimal and balanced presentation of the fundamental group of an orientable 3-dimensional small cover based on the Morse-theoretic approach. In section 5 we provide some illustrative examples including the dodecahedral space and the permutohedral space. We conclude with section 6 where the minimal Heegaard splitting is described. 

\section{Basics of small covers}

We recall some basic definitions and properties of small covers. For more information see \cite{DJ} or the monograph on toric topology \cite{BP}.

\begin{definition}
A small cover $M\xrightarrow{\pi}P$ over a simple $n$-dimensional convex polytope $P$ is a closed $n$-dimensional manifold $M$ with a locally standard $\mathbb{Z}_2^n$-action such that the orbit space $M/\mathbb{Z}_2^n$ is diffeomorphic to $P$ as a manifold with corners.
\end{definition}

Locally standard means that every point $p\in M$ has a stable neighborhood which is weakly equivariantly diffeomorphic with some stable open set in $\mathbb{R}^n$ with the standard action of $\mathbb{Z}_2^n$ generated by reflections in coordinate hyperplanes. The fixed points set of the action of $\mathbb{Z}^n_2$ on $M$ is the preimage $\pi^{-1}(V)$ of the set of vertices $V$ of the polytope $P$. The isotropy subgroup $G_F$ for a facet $F$ is a rank-one subgroup determined by a generator $\lambda(F)\in\mathbb{Z}_2^n$. In this way we obtain a map $\lambda:\mathcal{F}\rightarrow\mathbb{Z}_2^n$ on the set of facets $\mathcal{F}=\{F_1,\ldots,F_m\}$ of $P$, called the {\it characteristic map}. The pair $(P,\lambda)$ is the {\it characteristic pair} of the small cover $M\xrightarrow{\pi}P$. The isotropy subgroup $G_f$ for a codimension-$k$ face $f=F_{i_1}\cap\ldots\cap F_{i_k}$ is generated by $\lambda(F_{i_1}),\ldots,\lambda(F_{i_k})$.

As a quotient space the small cover $M$ is composed of $2^n$ copies of the polytope $P$ glued together along their facets
\begin{equation}\label{construction}
M=P\times\mathbb{Z}_2^n/\sim, \ \ \ \ \ (p,a)\sim(q,b) \ \mathrm{iff} \ p=q, \ b-a\in G_f,
\end{equation} where $G_f\subset\mathbb{Z}_2^n$ is the isotropy subgroup associated to a unique face $f\subset P$ such that $p$ lies in its relative interior $p\in f^{\circ}$.

The characteristic map $\lambda$ defines an $n\times m$-matrix $\Lambda$ over $\mathbb{Z}_2$ whose columns are the vectors $\lambda(F_i), i=1,\ldots,m$.
An arbitrary pair $(P,\lambda)$ consisting of a simple $n$-polytope $P$ and a map $\lambda:\mathcal{F}\rightarrow\mathbb{Z}_2^n$ is the characteristic pair of a small cover $M(P,\lambda)$ if and only if the following non-degeneracy condition holds
\[det \Lambda_v\neq 0,\] where $\Lambda_v$ is the minor of the characteristic matrix $\Lambda$ determined by all vectors $\lambda(F)$ with $v\in F$.

Small covers $M(P,\lambda)$ and $M(Q,\mu)$ are weakly equivariantly diffeomorphic if and only if there are a combinatorial equivalence of polytopes $\phi:P\rightarrow Q$ and a linear isomorphism $\theta:\mathbb{Z}_2^n\rightarrow\mathbb{Z}_2^n$ such that $\mu\phi=\theta\lambda$. Therefore the equivalence classes of characteristic pairs are in one-to-one correspondence with weak equivariant diffeomorphism classes of small covers. There are diffeomorphic small covers with non-equivalent characteristic pairs.

For a small cover $M=M(P,\lambda)$ each codimension-$k$ proper face $f\subset P$ defines the codimension-$k$ submanifold $M_f=\pi^{-1}(f)\subset M$ which is a small cover itself $M_f=M(f,\lambda_f)$, called the {\it facial submanifold}. Thereby the characteristic map $\lambda_f$ is defined in the following way. Let $\rho_f:\mathbb{Z}_2^n\rightarrow\mathbb{Z}_2^n/G_f\cong\mathbb{Z}_2^{n-k}$ be the projection onto the quotient group of $\mathbb{Z}_2^n$ by the isotropy subgroup $G_f$ of the face $f\subset P$. If $f=F_{i_1}\cap\ldots\cap F_{i_k}$ for some facets $F_{i_1},\ldots,F_{i_k}\in\mathcal{F}$ of $P$, then \[\lambda_f(F\cap f)=\rho_f(\lambda(F)),\] for each facet $F\in\mathcal{F}$ which intersects $f$ transversely, i.e. into a codimension one face of $f$. The non-degeneracy condition for $\lambda_f$ follows from that condition for $\lambda$.

\begin{example}[Example 1.20, \cite{DJ}]\label{surfaces}
The small covers over polygons are classified in the following way. For an $n$-gon $P_n$ a small cover $M=M(P_n,\lambda)$ is an orientable surface $S_g$ of genus $g$ if and only if $n=2g+2$ and $\lambda$ is a 2-valued characteristic map. In other cases the small cover $M$ over $P_n$ is a nonorientable surface $N_{n-2}$.
\end{example}

\begin{example}\label{simple}
The real projective space is a small cover over the simplex $\mathbb{R}P^n\xrightarrow{\pi}\Delta^n$ with the usual action of the multiplicative group $\mathbb{Z}_2^n$
\[(t_1,\ldots,t_n)[x_1:\ldots:x_n:x_{n+1}]=[t_1x_1:\ldots:t_nx_n:x_{n+1}],\] for $t_i\in\mathbb{Z}_2, i=1,\ldots,n, x_1^2+\ldots+x_{n+1}^2=1$. The characteristic matrix $\Lambda$ is $n\times(n+1)$-matrix with columns $\lambda(F_i)=e_i, i=1,\ldots,n$ and $\lambda(F_{n+1})=e_1+\cdots+e_n$, where $e_1,\ldots,e_n$ is the standard basis of $\mathbb{Z}_2^n$. The facial submanifold $M_f$ corresponding to a proper codimension-$k$ face $f\subset\Delta^n$ is equivariantly diffeomorphic to the real projective space $\mathbb{R}P^{n-k}$ with the usual action of $\mathbb{Z}_2^{n-k}$.

Similarly, the torus is a small cover over the cube $T^n\xrightarrow{\pi}I^n$
with the coordinate wise action of  $\mathbb{Z}_2^n$
\[(t_1,\ldots,t_n)(x_1,\ldots,x_n)=(t_1x_1,\ldots,t_nx_n), t_i\in\mathbb{Z}_2, x_i\in S^1, i=1,\ldots,n.\] The characteristic matrix is $n\times 2n$-matrix with columns $\lambda(F_i)=\lambda(F_{n+i})=e_i, i=1,\ldots,n$. The facial submanifolds are tori of smaller dimensions.
\end{example}

An orientability condition for small covers is given by Nakayama and Nishimura \cite{NN}. For a basis $\{e_1,\ldots,e_n\}$ of $\mathbb{Z}_2^n$ define a homomorphism $\epsilon:\mathbb{Z}_2^n\rightarrow\mathbb{Z}_2$ by $\epsilon(e_i)=1, i=1,\ldots,n$.

\begin{theorem}[Theorem 1.7, \cite{NN}]\label{orient} The small cover $M(P,\lambda)$ is orientable if and only if there is a basis $\{e_1,\ldots,e_n\}$ of $\mathbb{Z}_2^n$ such that the composition $\epsilon\lambda:\mathcal{F}\rightarrow\mathbb{Z}_2$ is the constant map $\epsilon(\lambda(F))=1, F\in\mathcal{F}$.
\end{theorem}

Special class of small covers, called {\it linear models} is given by $n$-colorable $n$-dimensional simple polytopes. For such polytopes a set of colors may be chosen to be a basis $\{e_1,\ldots,e_n\}$ of $\mathbb{Z}_2^n$. Any coloring $\lambda:\mathcal{F}\rightarrow\{e_1,\ldots,e_n\}$ gives rise to a characteristic map and defines the small cover $M=M(P,\lambda)$. Since any two colorings are equivalent up to reordering of colors the linear model is determined only by underlying polytope $M=M(P)$. Joswig proved the following combinatorial characterization of $n$-colorability of an $n$-polytope $P$.

\begin{theorem}[Theorem 16, \cite{J}]\label{even} A simple $n$-polytope $P$ is $n$-colorable if and only if all its 2-faces are evengons, the polygons with an even number of vertices.
\end{theorem}

The fundamental group of a small cover $M=M(P,\lambda)$ is determined in \cite[Lemma 4.5]{DJ} from the homotopy exact sequence of the Borel fibration

\[M\longrightarrow E\mathbb{Z}_2^3\times_{\mathbb{Z}_2^3}M\longrightarrow B\mathbb{Z}_2^3.\] The fundamnetal group $W_P=\pi_1(E\mathbb{Z}_2^3\times_{\mathbb{Z}_2^3}M)$ of the Borel construction of $M$ is the right-angled Coxeter group associated to the polytope $P$

\[W_P=\langle s_F,F\in\mathcal{F} | s_F^2=1, s_Fs_{F'}=s_{F'}s_F, F\cap F'\neq\emptyset\rangle.\] It determines the universal covering space $\mathcal{L}_P$ for all simple covers over the polytope $P$ by
\[\mathcal{L}_P=P\times W_P/\sim, \ (p,a)\sim (q,b) \ \mathrm{iff} \ p=q, ab^{-1}\in\langle s_F, p\in F\rangle.\] From the Borel fibration we obtain the short exaxt sequence
\[1\rightarrow\pi_1(M)\rightarrow W_P\xrightarrow{\phi}\mathbb{Z}_2^3\rightarrow 1,\] where $\phi$ is given by $\phi(s_F)=\lambda(F), \ F\in\mathcal{F}$. Therefore the fundamental group $\pi_1(M)$ is isomorphic to the kernel \[\pi_1(M)\cong\ker\phi.\]

A connected submanifold $X\subset M$ is called $\pi_1$-{\it injective} in a manifold $M$ if an embedding $i:X\subset M$ induces the monomorphism of fundamental groups $i_\ast:\pi_1(X)\rightarrow\pi_1(M)$. The following combinatorial criterium for $\pi_1$-injectivity is proved by Wu and Yu.

\begin{theorem}[Theorem 3.3 \cite{WY}]\label{injectivity}
The facial submanifold $M_f$ od a small cover $M=M(P,\lambda)$ is $\pi_1$-injective if and only if for any facet $F,F'\in\mathcal{F}$ transversal to a proper face $f\subset P$ we have $f\cap F\cap F'\neq\emptyset$ whenever $F\cap F'\neq\emptyset$.
\end{theorem}
As a consequence we obtain that the facet submanifold $M_F$ is $\pi_1$-injective if and only if a facet $F$ is not contained in any 3-{\it belt}, the set $\{F,F',F''\}\subset\mathcal{F}$ of three pairwise intersecting facets with the empty common intersection $F\cap F'\cap F''=\emptyset$. Recall that a simple polytope $P$ is called a {\it flag} polytope if any collection $F_{i_1},\ldots,F_{i_k}$ of pairwise intersecting facets $F_{i_r}\cap F_{i_s}\neq\emptyset, 1\leq r,s\leq k$ has the nonempty intersection $F_{i_1}\cap\ldots\cap F_{i_k}\neq\emptyset$. Therefore any face submanifold $M_f$ of a small cover $M=M(P,\lambda)$ over a flag polytope $P$ is $\pi_1$-injective. Note that this condition does not depend on the characteristic map $\lambda$, but only on the combinatorics of the polytope $P$.

Combining with the result of Davis, Januszkiewicz and Scott \cite{DJS} we obtain the following equivalent characterizations. Recall that a manifold is {\it aspherical} if its universal covering space is contractible.

\begin{theorem}[Theorem 2.2.5,\cite{DJS}, Proposition 3.6, \cite{WY}]\label{flag}
The following conditions are equivalent for a small cover $M=M(P,\lambda)$.
\begin{itemize}
\item $M$ is aspherical,
\item $P$ is a flag polytope,
\item All facial submanifolds $M_f\subset M$ are $\pi_1$-injective in $M$.
\end{itemize}
\end{theorem}

\subsection{3-dimensional small covers}

The 3-dimensional simple polytopes are represented by their Schlegel diagrams, which are characterized by the classical Steinitz theorem as 3-connected, 3-regular, planar graphs.

We are mainly concerned with orientable 3-dimensional small covers. The orientability condition from Theorem \ref{orient} in the dimension $n=3$ implies that
a 3-dimensional small cover $M=M(P,\lambda)$ is orientable if and only if the characteristic map $\lambda:\mathcal{F}\rightarrow\mathbb{Z}_2^3$ defines a 4-coloring of the polytope $P$ by colors $\{e_1,e_2,e_3,e_1+e_2+e_3\}$. Hence, by Four color theorem, every simple 3-polytope $P$ is the orbit space of an orientable 3-dimensional small cover $M\xrightarrow{\pi}P$.

The facet submanifold $M_F\subset M$ of a facet $F\subset P$ is the small cover $M_F=M(F,\lambda_F)$, where $\lambda_F:\{F'\in\mathcal{F} | F\cap F' \ \mathrm{is \ an \ edge \ of} \ P\}\rightarrow\mathbb{Z}_2^2$ is the characteristic map induced by the projection $\rho_F:\mathbb{Z}_2^3\rightarrow\mathbb{Z}_2^3/\langle\lambda(F)\rangle$
\[\lambda_F(F\cap F')=\rho_F(\lambda(F')).\] According to Example \ref{surfaces} the facet submanifold $M_F$ is an orientable surface of genus $g$ if and only if $F$ is a $(2g+2)$-gon and $\lambda_F$ is a 2-valued characteristic map. By Theorem \ref{even}, a polytope $P$ admits a linear model over it if and only if all its facets are evengons. Therefore for a linear model $M(P)$ each facet submanifold is an orientable surface.

The face submanifold $M_e\subset M$ is a circle consisting of two copies of the corresponding edge $e\subset P$. We denote these two copies by $e$ and $\overline{e}$, so if we orient the circle $M_e$ according to the edge $e$, we can write $M_e=e\overline{e}^{-1}$.

\begin{definition}
For a small cover $M=M(P,\lambda)$ let $M_F$ and $M_e$ be the face submanifolds corresponding to $2$-faces $F\subset P$ and to edges $e\subset P$. We call $M_F$ a {\it face surface} and $M_e$ an {\it edge circle} of $M$.
\end{definition}

The $\pi_1$-injectivity condition of Theorem \ref{injectivity} shows that the face surface $M_F$ of a 3-dimensional small cover $M=M(P,\lambda)$ is $\pi_1$-injective if and only if a face $F\subset P$ is not contained in any 3-belt. There always exists such a face, see Proposition 4.1 \cite{WY}.

Let us recall some basic definitions in the topology of 3-manifolds, for a detailed review see the classic monograph \cite{H}. We assume that all manifolds are closed and connected. Also for the sake of simplicity we assume that all manifolds are orientable. A manifold $M$ is called {\it prime} if it can not be represented as a connected sum $M=M_1\# M_2$ of two manifolds not diffeomorphic to the 3-sphere $S^3$. By the Milnor-Kneser theorem any manifold is uniquely decomposed into connected sum of prime manifolds up to insertion or deletion of $S^3$ summands. A manifold $M$ is called {\it irreducible} if every embedded sphere $S^2\subset M$ bounds a ball in $M$. The only prime manifold which is not irreducible is $S^1\times S^2$. A disc $D\subset M$ is called {\it essential} for an embedded, orientable, closed surface $S\subset M$, different from $S^2$, if $D$ meets $S$ transversally at the boundary $D\cap S=\partial D$ and $\partial D$ is a homotopy nontrivial loop on $S$. The surface $S\subset M$ is an {\it incompressible} in $M$ if it does not admit any essential discs. By the Loop theorem a surface $S$ is incompressible if and only if it is $\pi_1$-injective in $M$.

\begin{definition}
A closed, orientable, irreducible manifold $M$ is a {\it Haken manifold} if it contains a $\pi_1$-injective in $M$, orientable surface of positive genus.
\end{definition}

A well known and simple criteria for an irreducible manifold $M$ to be of the Haken type is given in the following lemma.

\begin{lemma}\label{betti}
A closed, orientable, irreducible manifold $M$ with the positive first Betti number $\beta_1(M)>0$ is a Haken manifold.
\end{lemma}

{\it Proof.} From the condition $\beta_1(M)>0$ we have that $H^1(M,\mathbb{Z})\cong H_1(M,\mathbb{Z})/Tors$ does not vanish. Since $H^1(M,\mathbb{Z})$ is represented by homotopy classes $[M,S^1]$, we conclude that there is a homotopy nontrivial map $f:M\rightarrow S^1$, which can be assumed to be smooth. Let $S=f^{-1}(\{x_0\})\subset M$ be an orientable surface defined by a regular value $x_0\in S^1$ of the map $f$. If there is an essential disc $D\subset M$ for $S$ we can cut the surface $S$ along this disc to obtain an orientable surface of smaller genus which represents the same nontrivial $2$-dimensional homology class. By iterating this procedure we end with a surface $S'$ of positive genus that is incompressible in $M$. The surface $S'$ can not be a sphere $S^2$ since $M$ is irreducible. \qed

The Haken manifolds play an important role in the topology of 3-manifolds. The existence of incompressible surfaces gives an inductive argument in proving the results about Haken manifolds. Waldhausen used this inductive procedure to prove that the fundamental group determines the diffeomorphism type of a Haken manifold.

\begin{theorem}[Corollary 6.5, \cite{W}]\label{fundamental}
Let $M$ and $N$ be closed, orientable, irreducible $3$-manifolds with isomorphic fundamental groups $\pi_1(M)\cong\pi_1(N)$. If additionally $M$ is a Haken manifold, then $M$ and $N$ are diffeomorphic.
\end{theorem}

We return to the class of small covers. The prime decomposition of an orientable small cover $M(P,\lambda)$ depends on the combinatorics of the polytope $P$.
A simple 3-polytope $P$ different from the simplex $\Delta^3$ is a flag polytope if and only if it has no 3-belts. Cutting along 3-belts decomposes an arbitrary simple polytope $P$ into a connected sum along vertices, whose summands are simplices and flag polytopes. This decomposition of the underlying polytope $P$ corresponds to the prime decomposition of the orientable small cover $M(P,\lambda)$ into an equivariant connected sum of real projective spaces $\mathbb{R}P^3$ and orientable small covers over flag polytopes. For details, see Theorem 4.3, \cite{E}.

\begin{example}
An orientable small cover over the 3-prism is equivariantly diffeomorphic to a connected sum $\mathbb{R}P^3\#\mathbb{R}P^3$. More generally, the truncation of a vertex of a 4-colored simple polytope $P$ corresponds to an equivariant connected sum $M\#\mathbb{R}P^3$ of an orientable small cover $M$ over $P$. Since truncations of different vertices do not have to lead to combinatorially equivalent polytopes, we obtain a series of examples of diffeomorphic but not equivariantly diffeomorphic small covers.
\end{example}

\begin{remark}\label{pogorelov}
A special class of flag polytopes are Pogorelov polytopes which can be realized in Lobachevsky space as right-angled polytopes. They are characterized in works of Pogorelov and Andreev as flag polytopes without 4-belts. The small covers over Pogorelov polytopes are hyperbolic 3-manifolds.
\end{remark}

We collect the topological properties of small covers in the following proposition. This was proved regardless to orientability in the unpublished arXiv version \cite{WYa} of the paper \cite{WY}.
\begin{proposition}\label{haken}
An orientable small cover $M\xrightarrow{\pi}P$ over a flag polytope $P$ is a Haken manifold.
\end{proposition}
{\it Proof.} The small cover $M$ is an irreducible manifold since $P$ is a flag polytope. We need to find a $\pi_1$-injective in $M$, orientable surface. If there is an orientable face surface $M_F$, we done since it is $\pi_1$-injective by Theorem \ref{flag}. If there are not, we choose a facet $F$ and its small neighborhood $U_F$ in $P$. The inverse image $\pi^{-1}(U_F)$ under the orbit map is a tubular neighborhood of the face surface $M_F$. Since $M_F$ is one-sided in $M$, the boundary of this tubular neighborhood is an orientable surface $M_{\widetilde{F}}$ which double covers the surface $M_F$. Let $\eta:M_{\widetilde{F}}\rightarrow M_F$ be the double covering and $i_1:M_F\rightarrow M$ and $i_2:M_{\widetilde{F}}\rightarrow M$ be embeddings. Then $i_1\eta$ induces a monomorphism of fundamental groups and the composition $i_1\eta$ is homotopic to $i_2$ via a deformation retraction of the tubular neighborhood onto the surface $M_F$. Thus the orientable surface $M_{\widetilde{F}}$ is a $\pi_1$-injective in $M$. \qed

By Waldhausen Theorem \ref{fundamental} and Proposition \ref{haken} we obtain that the fundamental groups classify orientable small covers over flag polytopes. Since every factor of the prime decomposition of an orientable small cover $M$ is $\mathbb{R}P^3$ and small covers over flag polytopes, the following theorem is a consequence of a unique factorization of finitely presented groups into indecomposable factors according to the free product of groups.

\begin{theorem}[Corollary 4.10, \cite{WYa}]
The orientable 3-dimensional small covers $M_1$ and $M_2$ with the isomorphic fundamental groups $\pi_1(M_1)\cong\pi_1(M_2)$ are diffeomorphic.
\end{theorem}

\section{The canonical Heegaard splitting}

A convenient way to compute the fundamental group of a closed 3-manifolds is obtained from Heegaard splittings. A {\it handlebody} $H$ of genus $g$ is a connected, orientable 3-manifold with boundary containing a collection of properly embedded 2-discs $D_1,\ldots, D_g\subset H$, called {\it essencial discs}, such that cutting along these discs produces a 3-ball $H\setminus\bigcup D_i=B^3$. The boundary $\partial H$ is an orientable surface of genus $g$. A {\it Heegaard splitting} $(\Sigma,H,H')$ of genus $g$ of a closed, orientable 3-manifold $M$ is a decomposition into two handlebodies $M=H\cup H'$ of genus $g$ such that $\Sigma=H\cap H'$ is their common boundary surface $\partial H=\partial H'=\Sigma$. The {\it Heegaard genus} $\gamma(M)$ is the minimal genus of Heegaard splittings of $M$.

Given a Heegaard splitting $(\Sigma,H,H')$ of genus $g$ of a manifold $M$, denote by $\boldsymbol{\mu}=\{\mu_1,\ldots,\mu_g\}$ and $\boldsymbol{\lambda}=\{\lambda_1,\ldots,\lambda_g\}$ the systems of simple, closed curves on the surface $\Sigma$ which are intersections of $\Sigma$ with collections of essential discs in the handlebodies $H$ and $H'$ respectively. The {\it Heegaard diagram} is the triplet $(\Sigma,\boldsymbol{\mu},\boldsymbol{\lambda})$ assuming that the curves from the corresponding systems intersect transversally. The systems of curves $\boldsymbol{\mu},\boldsymbol{\lambda}$ we call meridian-longitude systems. Cutting the surface $\Sigma$ along the meridian system $\boldsymbol{\mu}$ produces a 2-sphere with paired $2g$ holes. Then the longitude system $\boldsymbol{\lambda}$ transforms into a system of arcs connecting these holes. The manifold $M$ can be reconstructed from these data. First we glue paired holes on the sphere, seen as the boundary of the ball $B^3$, and obtain a genus $g$ handlebody $H$ with meridian discs. The arcs, previously obtained from the system of curves $\boldsymbol{\lambda}$, connect into closed, simple curves on the boundary $\partial H$. Gluing 2-handles along these curves and afterthat adding the remaining 3-ball reconstruct the manifold $M$.

Giving a Heegaard diagram of a manifold $M$ the fundamental group $\pi_1(M)$ can be presented by the Van Kampen theorem with

\begin{equation}\label{Heg}
\pi_1(M)=\langle s_1,\ldots,s_g | r_1,\ldots,r_g\rangle,
\end{equation} where $s_1,\ldots,s_g$ are the generators corresponding to the cores of handles of $H$ and $r_1,\ldots,r_g$ are relations obtained in the following way. Suppose the curves from the meridian-longitude systems are oriented. We go in the direction of a longitude curve $\lambda_i$ and each time when pass across a meridian curve $\mu_j$ associate a letter $s_j^{\pm 1}$ according to the local intersection number. The word we obtain gives the relation $r_i$. It follows that the fundamental groups of oriented 3-manifolds are balanced, i.e. allow presentations with the same numbers of generators and relations. The rank $r(G)$ of a finitely presented group $G$ is the minimal number of generators needed to present $G$. We obtain
\begin{equation}\label{rank}
r(\pi_1(M))\leq\gamma(M),
\end{equation} and the equality is not achieved in general case.

We construct the canonical Heegaard diagram of an orientable small cover $M=M(P,\lambda)$ out off the characteristic pair $(P,\lambda)$. Denote by $P_g, g\in\mathbb{Z}_2^3$ the copies of the polytope $P$ and by $F_g, F\in\mathcal{F}, g\in\mathbb{Z}_2^3$ the collections of their facets. According to the construction $(\ref{construction})$ of $M(P,\lambda)$, we glue these copies along the facets by the rules
\[F_g\sim F_{g+\lambda(F)}, F\in\mathcal{F}, g\in\mathbb{Z}_2^3.\] Let $U$ be the thickened $1$-skeleton of the polytope $P$ and $V$ its complement $V=P\setminus U$. Performing the gluing process on their copies $U_g, V_g, g\in\mathbb{Z}_2^3$ as indicated by the coloring $\lambda$ produces the handlebodies $H_1=V\times\mathbb{Z}_2^3/\sim$ and $H_2=U\times\mathbb{Z}_2^3/\sim$ which gives a Heegaard splitting of $(\Sigma,H_1,H_2)$ of the small cover $M(P,\lambda)$. We say that this splitting is a {\it canonical} for $M(P,\lambda)$. Denote by $N^{n,m}$ a handlebody formed from $n$ balls $B$ and $m$ handles $I\times D^2$. Then $H_1=N^{8,4f_2}$ and $H_2=N^{f_0,2f_1}$, where $f_0,f_1,f_2$ are the numbers of vertices, edges and faces of the polytope $P$. Note that Euler characteristics of $H_1$ and $H_2$ are equal $8-4f_2=f_0-2f_1$ by the Dehn-Sommerville relation for $P$.

The handlebody $H_1$ is formed from eight balls $B_g$ lying in the interiors of $P_g$ and from handles which correspond to the pairs of glued facets $(F_g,F_{g+\lambda(F)}), F\subset P, g\in\mathbb{Z}_2^3$. Denote by $s_{F,g}$ the core of the handle corresponding to $(F_g,F_{g+\lambda(F)})$. Note that $s_{F,g}$ and $s_{F,g+\lambda(F)}$ are the same curves with the opposite orientations, so we can write $s_{F,g+\lambda(F)}=s_{F,g}^{-1}$. Similarly, the handlebody $H_2$ is formed from balls $B_v, v\in V$, where $V$ is the vertex set of $P$ and from handles, one handle for each pair consisting of an edge $e$ and a coset in $\mathbb{Z}_2^3/G_e$. A handle $(e,[g])$ corresponding to an edge $e=F\cap F'$ and a coset $[g]\in\mathbb{Z}_2^3/\langle\lambda(F),\lambda(F')\rangle$ passes across handles of $H_1$ giving the word $s_{F,g}s_{F',g+\lambda(F)}s_{F,g+\lambda(F')}^{-1}s_{F',g}^{-1}$ that is determined by a 4-cycle of glued faces
\[F_g\sim F_{g+\lambda(F)}, F'_{g+\lambda(F)}\sim F'_{g+\lambda(F)+\lambda(F')},\] \[F_{g+\lambda(F)+\lambda(F')}\sim F_{g+\lambda(F')}, F'_{g+\lambda(F')}\sim F'_g.\]

Choose a vertex $v_0\in P$. The eight balls $B_g,g\in\mathbb{Z}_2^3$ and handles with cores $s_{F,g}, v_0\in F$ form the configuration of the $1$-skeleton of the 3-cube. Adding the ball $B_{v_0}$ and handles from $H_2$ corresponding to edges that contain the vertex $v_0$ gives a neighborhood $D$ of $v_0$ in $M$ homeomorphic to a ball. The ball $D$ and remaining handles of $H_1$ form a new handlebody $H'_1\subset M$ homeomorphic to $N^{1,4(f_2-3)}$. The complement $H'_2=M\setminus H'_1$ is a handlebody $H_2$ with removed ball $B_{v_0}$ and handles corresponding to edges that contain the vertex $v_0$, thus $H'_2$ is homeomorphic to $N^{f_0-1,2f_1-6}$. Finally we find a maximal tree $T$ of the $1$-skeleton of $P$ with the vertex $v_0$ being a leaf. The balls $B_v,v\neq v_0$ and handles corresponding to edges of the tree $T\setminus v_0$ form a ball $D'$ in $H'_2$. The ball $D'$ and remaining handles of $H'_2$ form the handlebody $H''_2$ homeomorphic to $N^{1,2f_1-6-(f_0-2)}$. We obtain a new Heegaard splitting $M=H'_1\cup H''_2$ of the small cover $M$ whose genus is equal to $4(f_2-3)$. From the canonical Heegaard splitting we rediscover the Wu and Yu presentation of the fundamental group of the orientable 3-dimensional small cover $M$.
\begin{proposition}[Proposition 2.1, \cite{WY}]\label{present}
The fundamental group $\pi_1(M(P,\lambda),v_0)$ is generated by $s_{F,g}, F\in\mathcal{F}, g\in\mathbb{Z}_2^3$ with the following set of relations

\[s_{F,g+\lambda(F)}=s_{F,g}^{-1}, F\in\mathcal{F}, g\in\mathbb{Z}_2^3,\]
\[s_{F,g}s_{F',g+\lambda(F)}=s_{F',g}s_{F,g+\lambda(F')}, F\cap F'\neq\emptyset,\]
\[s_{F,g}=1, v_0\in F.\]
\end{proposition}

Note that $f_0-2$ of the middle relations in Proposition \ref{present} are redundant after choosing a maximal tree $T$, so the presentation can be make to be  balanced.

\section{The minimal presentation of $\pi_1(M)$}

The presentation of the fundamental group $\pi_1(M)$ given by Proposition \ref{present} have the rank $4(f_2-3)$. Since there is an epimorphism \[\pi_1(M)\rightarrow H_1(M,\mathbb{Z})\rightarrow H_1(M,\mathbb{Z}_2),\] we have that the rank of $\pi_1(M)$ is bounded from bellow by the first mod 2 Betti number $r(\pi_1(M))\geq\beta^{\mathbb{Z}_2}_1(M)$. The mod 2 Betti numbers of small covers are equal to the components of the $h$-vector of $P$ \cite[Theorem 3.1]{DJ}, so $\beta^{\mathbb{Z}_2}_1(M)=f_2-3$. We construct a Heegaard splitting of $M$ of genus $f_2-3$ and use it to present the fundamental group $\pi_1(M)$. It follows from $(\ref{rank})$ that
\[r(\pi_1(M))=\gamma(M)=f_2-3,\] which shows that the obtained presentation will be minimal.

The following construction, based on the Morse-theoretic approach, is standard and firstly used by Khovanskii \cite{K} for determining the canonical cell structure of smooth toric varieties, for more details see \cite{DJ}, \cite{BP}. Given a small cover $M=M(P,\lambda)$ over a simple $n$-dimensional polytope $P\subset\mathbb{R}^n$, let $\nu\in\mathbb{R}^n$ be a generic vector according to $P$. This means that $\langle \nu,v_i-v_j\rangle\neq 0$ for each pair of vertices $v_i,v_j\in V$ of the polytope $P$, where $\langle,\rangle$ is an inner product on $\mathbb{R}^n$. Then the function $\varphi(x)=\langle \nu,\pi(x)\rangle, x\in M$, can be modified into a Morse function $f:M\rightarrow\mathbb{R}$ on the small cover $M\xrightarrow{\pi}P$. We first smoothly modify the functional $\langle \nu,\cdot\rangle$ to obtain a vector field on $P$ which is tangent to each face and vanish in vertices of $P$. The gradient vector field of $f$ is the pull-back of such vector field on $P$ under the projection map $\pi:M\rightarrow P$. The critical points of the function $f$ are the inverse images $\overline{v}=\pi^{-1}(v)\in M, v\in V$ of the vertices of $P$.

The functional $\langle \nu,\cdot\rangle$ can be seen as a height function on $P$. The orienting each edge into a direction of increasing height turns the $1$-skeleton $P^{(1)}$ into an oriented graph. The index $\mathrm{ind}_\nu(v)$ of a vertex $v\in P$ is defined as the number of incoming edges to $v$. Each vertex is maximal for a unique face $G_v\subset P$, spanned by incoming edges of $v$ and minimal for a unique face $G^v\subset P$, spanned by outgoing edges. The number of vertices of index $i$ is equal to the $i$-th component $h_i$ of the $h$-vector of $P$. The Morse index of a critical point $\overline{v}\in M$ is equal to the index $\mathrm{ind}_\nu(v)$ of the vertex $v\in V$. The facial submanifolds $M_{G_v}, M_{G^v}\subset M$ are closures of the descending and ascending manifolds of the critical point $\overline{v}\in M$. Note that, generally, the obtained Morse function is not of the Morse-Smale type.

Descending submanifolds gives a cellular subdivision of $M$ which is perfect with respect to $\mathbb{Z}_2$ coefficients. The fundamental group $\pi_1(M)$ depends only on the $2$-skeleton of this subdivision. Since the critical points of the function $f$ in general are not arranged in increasing order of their indices, it could be hard to control the behavior of flow-lines in dimensions higher than $3$. In dimension $3$ we can handle it which is the basic for the following theorem.

\begin{theorem}\label{main}
Let $M\xrightarrow{\pi}P$ be an oriented $3$-dimensional small cover over a simple polytope $P\subset\mathbb{R}^3$ and $h(x)=\langle\nu,x\rangle, x\in P$ be a height function for a vector $\nu\in\mathbb{R}^3$ generic according to $P$. Let $v_0\in P$ be the vertex with the minimal height. Then the fundamental group $\pi_1(M,\overline{v}_0)$ in the base point $\overline{v}_0=\pi^{-1}(v_0)$ has a minimal presentation with the generators corresponding to vertices of index $1$ and the relations corresponding to vertices of index $2$.
\end{theorem}

{\it Proof}. A Morse function $f:M\rightarrow\mathbb{R}$, constructed in a way described above from the height function $\langle\nu,\cdot\rangle$ on $P$, determines a handle decomposition of $M$. The vertex $v_0$ is the source of the $1$-skeleton of $P$ whose edges are oriented in directions of increasing the height function. There are the same number of $f_2-3$ vertices of indices $1$ and $2$, where $f_2$ is the number of faces of $P$. We label edges and faces excluding those that contain the sink.

Choose a directly oriented simple edge path $p_v=e_0\ldots e_v$ from the source $v_0$ to each vertex $v$ of index one. The correspondence of a vertex $v$ of index one to the final edge $e_v$ of the edge path $p_v$ is unambiguously determined. The edge loop $\alpha_v=(p_v\setminus e_v)M_{e_v}(p_v\setminus e_v)^{-1}$, where $M_{e_v}$ is the edge circle, represents a generator of the fundamental group of $M$ at the base point $v_0$.

The face surface $M_F$ is a small cover $M_F=M(F,\lambda_F)$ built off four copies of a face $F\subset P$ glued along edges by the rule determined by the coloring $\lambda_F$. To each face $F\subset P$ is associated a unique vertex $v_F$ with the maximal height. We choose a shelling order on the set of faces starting with those containing the source and ending with those containing the sink. In this order, each face $F\subset P$, excluding the last three, determines the relation $r_F$. The relation $r_F$ is read off the boundary of the $2$-cell obtained by gluing four copies of $F$ along edges that contain the vertex $v_F$. \qed

\begin{remark}
The procedure for calculation of the fundamental group of an orientable $3$-dimensional small cover depends on the choices of a height function and of edge paths
from the source vertex to vertices of index one, so we obtain different presentation depending on these choices. The same procedure is valid in the case of nonorientable small cover.
\end{remark}

\section{Examples}

We illustrate the procedure described by Theorem \ref{main} in several examples.

\begin{example}[The projective space $\mathbb{R}P^3$]
The most simple example is the small cover $\mathbb{R}P^3\xrightarrow{\pi}\Delta^3$, see Example \ref{simple}. The oriented edge graph defined by a height function on the simplex is shown in Figure \ref{proj}. The simplex is $4$-colored, we label its edges and faces not containing the sink. The edge circle $M_a$ generates $\pi_1(\mathbb{R}P^3,v_0)$ with the only relation corresponding to the face $F$ numerated by $1$, whose face surface $M_F$ is the projective plane $\mathbb{R}P^2$.
\begin{figure}[h!h!h!]
\begin{center}
\includegraphics[width=100mm,height=60mm]{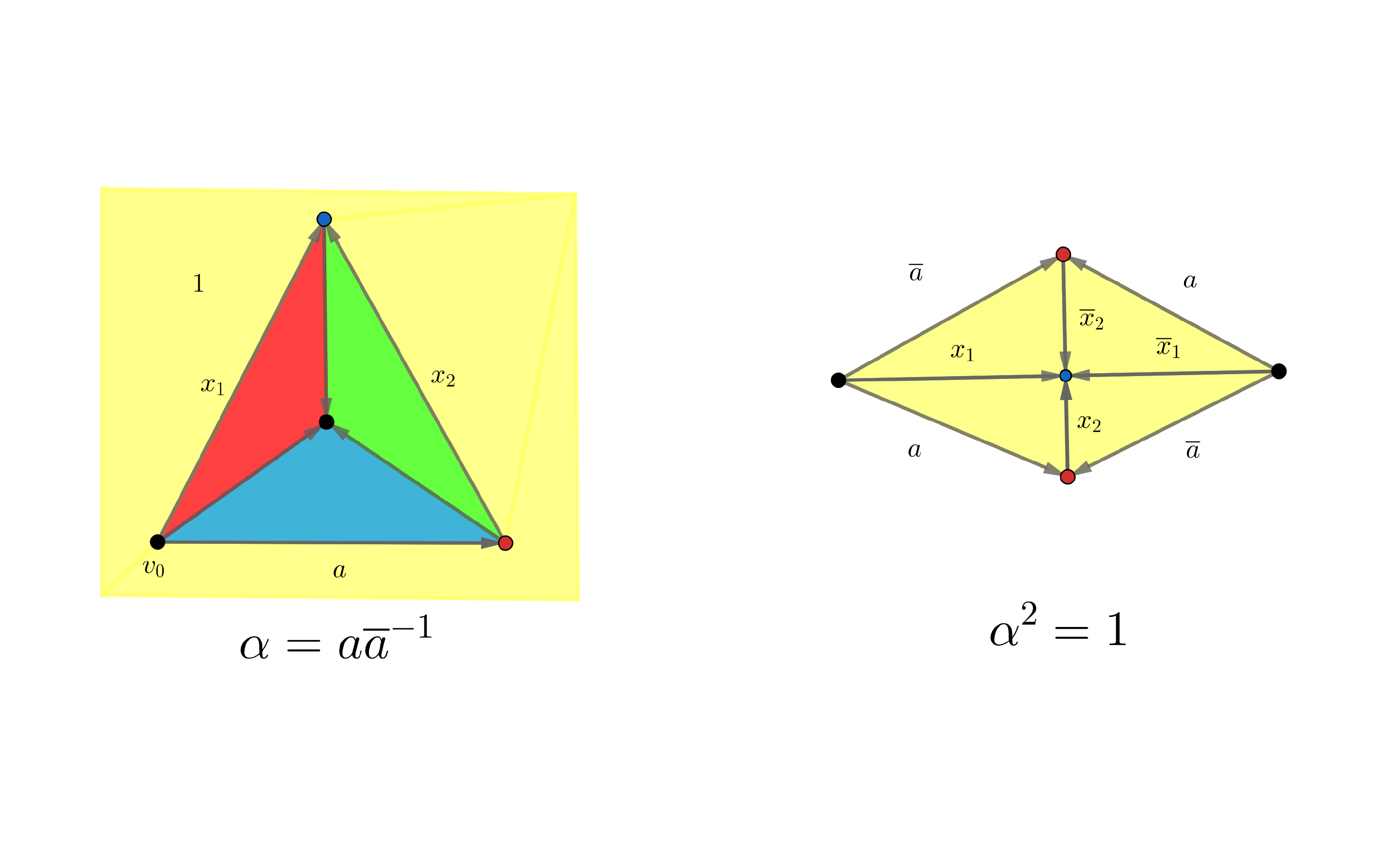}
\caption{Projective space}\label{proj}
\end{center}
\end{figure}
\end{example}

\begin{example}[Dodecahedral space]

In Figure \ref{dodec} is given a coloring of the dodecahedron $D$. Since $D$ is a Pogorelov polytope, this produces an orientable small cover $M\xrightarrow{\pi}D$ which is a Haken, hyperbolic manifold.
\begin{figure}[h!h!h!]
\begin{center}
\includegraphics[width=160mm]{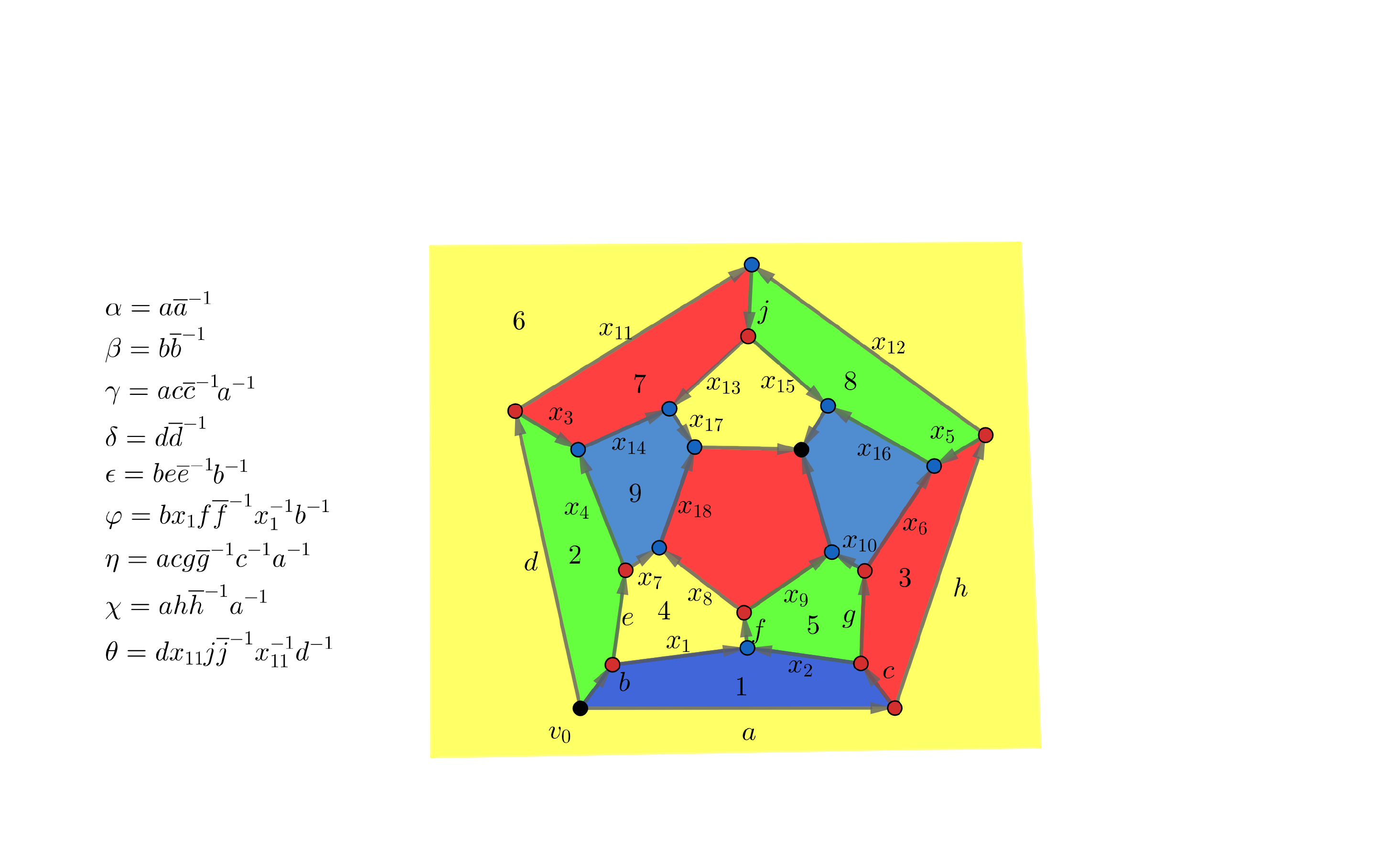}
\caption{Dodecahedral space}\label{dodec}
\end{center}
\end{figure}

We orient edges according to a height function and add them labels as indicated. The vertices of index $1$ are colored in red and vertices of index $2$ are colored in blue. The edge loops that generate the fundamental group $\pi_1(M,v_0)$ correspond to edges ending in red vertices. The relations are given by numerated faces and each of them has the maximal height in a unique blue vertex. Note that each face surface is the nonorientable $N_3=T^2\#\mathbb{R}P^2$. For the sake of certainty, suppose that colors correspond to the vectors of $\mathbb{Z}_2^3$ in the following way: blue $e_1$, green $e_2$, yellow $e_3$ and red $e_1+e_2+e_3$. The reading the relation out of a numerated face is illustrated by the faces $1$ and $2$ in Figure \ref{relD}.

Calculation gives the following relations in generators defined on Figure \ref{dodec}. The relations are listed according to numeration of faces:
\[\gamma\alpha\beta\gamma^{-1}\alpha\beta^{-1}=1, \epsilon\beta\delta\epsilon^{-1}\beta\delta^{-1}=1, \eta\gamma\chi\gamma^{-1}\eta\chi^{-1},\]
\[\epsilon\alpha^{-1}\gamma^{-1}\varphi^{-1}\gamma\alpha\epsilon\varphi^{-1}=1, \eta\alpha\beta^{-1}\varphi^{-1}\beta\alpha^{-1}\eta\varphi^{-1}=1,\]
\[\chi\delta^{-1}\alpha^{-1}\chi^{-1}\alpha\delta^{-1}=1, \theta^{-1}\chi\delta^{-1}\beta^{-1}\epsilon^{-1}\delta\chi^{-1}\theta^{-1}\epsilon\beta=1,\]
\[\alpha\delta^{-1}\theta\gamma^{-1}\eta^{-1}\theta\delta\alpha^{-1}\gamma^{-1}\eta^{-1}=1, \varphi\epsilon^{-1}\beta\chi^{-1}\theta^{-1}\delta\beta^{-1}\varphi\epsilon^{-1}\delta\chi^{-1}\theta^{-1}=1.\]

Abelianization of the obtained group gives $H_1(M,\mathbb{Z})\cong\mathbb{Z}_2^9$. From the Suciu-Trevisan formula \cite{ST} for rational Betti numbers we get $\beta_1(M)=0$ which means that the dodecahedral space $M$ is a rational homology sphere. Therefore the example of the dodecahedral space $M$ shows that the condition $\beta_1(M)>0$, stated in Lemma \ref{betti}, is not necessary for a manifold to be a Haken.

\end{example}

\begin{example}[Permutohedral space]

The permutohedron $P$ is $3$-colorable simple polytope by Theorem \ref{even}, so it determines a unique linear model $M(P)\xrightarrow{\pi}P$. We derive a presentation of the fundamental group $\pi_1(M(P),v_0)$ from the Bruhat order on the permutation group $S_4$. For the set of generators  $\{A,B,C,D,E,F,G,H,I,J,K\}$ we obtain the following relations:
\[BC=CB, DE=ED, DB^{-1}FBA^{-1}D^{-1}AF^{-1}=1, GA^{-1}EAC^{-1}G^{-1}CE^{-1}=1, \]
\[HDB^{-1}H^{-1}BD^{-1}=1, IEC^{-1}I^{-1}CE^{-1}=1, JAF^{-1}HFA^{-1}E^{-1}J^{-1}EH^{-1}=1,\]
\[JAG^{-1}IGA^{-1}D^{-1}J^{-1}DI^{-1}=1, HC^{-1}KCF^{-1}H^{-1}FK^{-1}=1,\]
\[KG^{-1}IGB^{-1}K^{-1}BI^{-1}=1, HC^{-1}IB^{-1}CH^{-1}BI^{-1}=1.\]

Abelianization gives $H_1(M(P),\mathbb{Z})\cong\mathbb{Z}^{11}$, which also follows from the fact that all facial submanifolds of $M(P)$ are orientable, see \cite[Corollary 3.7]{DJ}.
\end{example}

\section{The minimal Heegaard diagram}

The presentation $(\ref{Heg})$ of the fundamental group of a closed, oriented $3$-manifold $M$ depends on the Heegaard splitting. In return, having a balanced presentation of the fundamental group of $M$, we can construct a system of simple, closed curves on the orientable surface of genus that is equal to the rank of the presentation, whose curves correspond to relations. This system of curves defines a Heegaard diagram for some closed, orientable $3$-manifold $N$ with the fundamental group isomorphic to that of $M$, but nevertheless which may not be diffeomorphic to $M$. The counterexample is given by homotopic but not homeomorphic lens spaces. However, the Waldhausen theorem \ref{fundamental} shows that in the case of $M$ being a Haken manifold, the obtained system of curves is in fact a Heegaard diagram of $M$. Consequently, we can use Theorem \ref{main} to construct a Heegaard diagram of the smallest genus of an orientable small cover $M=M(P,\lambda)$.





\begin{thebibliography}{9}

\bibitem{BEMPP}  V. Buchstaber, N.\,Erokhovets, M. Masuda, T. Panov, S. Park, \textit{Cohomological rigidity of manifolds defined by 3-dimensional polytopes}, Russian Mathematical Surveys {\bf 72} (2017), no. 2, 199--256.

\bibitem{BP} V. Buchstaber, T. Panov, \textit{Toric Topology}, Mathematical Survays and Monographs, 204. Amer. Math. Soc., Providence, RI, (2015)

\bibitem{DJ} M. Davis, T. Januszkiewitc, {\it Convex polytopes, Coxeter orbifolds and torus actions}, Duke Math. J. {\bf 62} (1991), no.2, 417--451.

\bibitem{DJS} M. Davis, T. Januszkiewitc, R. Scott, \textit{Nonpositive curvature of blow-ups}, Selecta Math. \textbf{4} (1998), 491--547.

\bibitem{E} N. Erokhovets, \textit{Canonical geometrization of 3-manifolds realizable as small covers}, {arXiv:2011.11628}

\bibitem{H} J. Hempel, \textit{3-Manifolds}, Annals of Math. Studies, vol. 86, Princeton University Press, (1976)

\bibitem{J} M. Joswig, \textit{Projectivities in simplicial complexes and colorings of simple polytopes}, Math. Z. \textbf{240} (2002), 243--259.

\bibitem{K} A. Khovanskii, \textit{Hyperplane sections of polyhedra, toric varieties and discrete groups in Lobachevskii space}, Funk. Anal. i Prilozhen. \textbf{20} (1986), no.1, 50--61. (Russian); Funct. Anal. Appl. \textbf{20} (1986), no.1, 41--50. (English translation)

\bibitem{NN} H. Nakayama, Y. Nishimura, \textit{The orientability of small covers and coloring simple polytopes}, Osaka J. Math. \textbf{42} (2005), 243--256.

\bibitem{ST} A. Suciu, A. Trevisan, \textit{Real toric varieties and abelian covers of generalized Davis-Januszkiewicz spaces}, preprint 2012.

\bibitem{WY} L. Wu, L. Yu, \textit{Fundamental groups of small covers revisited}, Int. Math. Res. Not., {\bf 10} (2021), 7262--7298.

\bibitem{WYa} L. Wu, L. Yu, \textit{Fundamental groups of small covers revisited}, {arXiv:1712.00698v2}

\bibitem{W} F. Waldhausen, \textit{On irreducible 3-manifolds which are sufficiently large}, Ann. Math. {\bf 87} (1968), 56--88.

\end{thebibliography}
\end{document}